\documentclass[12pt]{article}
\usepackage{amssymb,latexsym,pstcol,pst-plot,amsthm,amsmath,hyperref}
\usepackage[hmargin=1in,vmargin=1in]{geometry}

\newcommand{\bal}{\begin{align*}}
\newcommand{\eal}{\end{align*}}
\newcommand{\gauss}[2]{\genfrac{[}{]}{0pt}{}{#1}{#2}} 

\newcommand{\Gaaea}{\put(0,0){\line(1,0){40}}}

\newcommand{\Gaaaf}{\put(0,0){\line(0,1){50}}}

\newcommand{\Gbabf}{\put(10,0){\line(0,1){50}}}

\newcommand{\Gcacf}{\put(20,0){\line(0,1){50}}}

\newcommand{\Gdadf}{\put(30,0){\line(0,1){50}}}

\newcommand{\Geaef}{\put(40,0){\line(0,1){50}}}

\newcommand{\Gabeb}{\put(0,10){\line(1,0){40}}}

\newcommand{\Gacec}{\put(0,20){\line(1,0){40}}}

\newcommand{\Gaded}{\put(0,30){\line(1,0){40}}}

\newcommand{\Gaeee}{\put(0,40){\line(1,0){40}}}

\newcommand{\Gafef}{\put(0,50){\line(1,0){40}}}

\newcommand{\ben}{\begin{enumerate}}
\newcommand{\een}{\end{enumerate}}
\newcommand{\ble}{\begin{lem}}
\newcommand{\ele}{\end{lem}}
\newcommand{\bth}{\begin{thm}}
\renewcommand{\eth}{\end{thm}}
\newcommand{\bpr}{\begin{prop}}
\newcommand{\epr}{\end{prop}}
\newcommand{\bco}{\begin{cor}}
\newcommand{\eco}{\end{cor}}
\newcommand{\bcon}{\begin{conj}}
\newcommand{\econ}{\end{conj}}
\newcommand{\bde}{\begin{defn}}
\newcommand{\ede}{\end{defn}}
\newcommand{\bex}{\begin{exa}}
\newcommand{\eex}{\end{exa}}
\newcommand{\barr}{\begin{array}}
\newcommand{\earr}{\end{array}}
\newcommand{\btab}{\begin{tabular}}
\newcommand{\etab}{\end{tabular}}
\newcommand{\beq}{\begin{equation}}
\newcommand{\eeq}{\end{equation}}
\newcommand{\bea}{\begin{eqnarray*}}
\newcommand{\eea}{\end{eqnarray*}}
\newcommand{\bce}{\begin{center}}
\newcommand{\ece}{\end{center}}
\newcommand{\bpi}{\begin{picture}}
\newcommand{\epi}{\end{picture}}
\newcommand{\bpp}{\begin{picture}}
\newcommand{\epp}{\end{picture}}
\newcommand{\bfi}{\begin{figure} \begin{center}}
\newcommand{\efi}{\end{center} \end{figure}}
\newcommand{\bprf}{\begin{proof}}
\newcommand{\eprf}{\end{proof}\medskip}
\newcommand{\capt}{\caption}
\newcommand{\bsl}{\begin{slide}{}}
\newcommand{\esl}{\end{slide}}
\newcommand{\bfr}{\begin{frame}}
\newcommand{\efr}{\end{frame}}

\newcommand{\prf}{{\bf Proof.\hspace{7pt}}}
\newcommand{\hqed}{\hfill \qed}
\newcommand{\hqedm}{\hfill \qed \medskip}

\newcommand{\eqed}[1]{$\textcolor{white}{\qed}\hfill{\dil#1}\hfill\qed$}

\newcommand{\ol}{\overline}

\newcommand{\hs}[1]{\hspace{#1}}
\newcommand{\hso}[1]{\hspace{-1pt}}
\newcommand{\vs}[1]{\vspace{#1}}
\newcommand{\qmq}[1]{\quad\mbox{#1}\quad}

\newcommand{\pmo}[2]{\put(#1){\makebox(0,0){#2}}}

\newcommand{\sbs}{\subset}
\newcommand{\sbe}{\subseteq}

\newcommand{\fl}[1]{\lfloor #1 \rfloor}

\def\<{\langle}
\def\>{\rangle}

\newcommand{\ree}[1]{(\ref{#1})}

\newcommand{\ra}{\rightarrow}

\newcommand{\al}{\alpha}
\newcommand{\be}{\beta}
\newcommand{\ga}{\gamma}
\newcommand{\de}{\delta}
\newcommand{\ep}{\epsilon}
\newcommand{\vep}{\varepsilon}

\newcommand{\ka}{\kappa}
\newcommand{\la}{\lambda}

\newcommand{\om}{\omega}

\newcommand{\si}{\sigma}

\newcommand{\La}{\Lambda}

\newcommand{\bbN}{{\mathbb N}}
\newcommand{\bbP}{{\mathbb P}}

\newcommand{\bbZ}{{\mathbb Z}}

\newcommand{\cB}{{\cal B}}

\newcommand{\cF}{{\cal F}}
\newcommand{\cG}{{\cal G}}
\newcommand{\cH}{{\cal H}}

\newcommand{\fS}{{\mathfrak S}}

\newcommand{\Av}{\operatorname{Av}}

\newcommand{\des}{\operatorname{des}}
\newcommand{\Des}{\operatorname{Des}}

\newcommand{\inv}{\operatorname{inv}}
\newcommand{\Inv}{\operatorname{Inv}}

\newcommand{\maj}{\operatorname{maj}}

\newcommand{\sta}{\operatorname{st}}

\newcommand{\bul}{\bullet}
\newcommand{\dil}{\displaystyle}

\newcommand{\x}{\hs{-2pt}\bul\hs{-2pt}}
\newcommand{\gab}{\ol{\ga}}
\newcommand{\bin}[2]{\left\{ #1 \atop #2 \right\}}

\newcommand{\exc}{\mathop{\rm exc}}

\newtheorem{thm}{Theorem}[section]
\newtheorem{prop}[thm]{Proposition}
\newtheorem{cor}[thm]{Corollary}
\newtheorem{lem}[thm]{Lemma}
\newtheorem{conj}[thm]{Conjecture}
\newtheorem{exa}[thm]{Example}
\newtheorem{question}[thm]{Question}

\begin{document}
\pagestyle{plain}

\title{Mahonian pairs
}
\author{
Bruce E. Sagan\\[-5pt]
\small Department of Mathematics, Michigan State University,\\[-5pt]
\small East Lansing, MI 48824-1027, USA, \texttt{sagan@math.msu.edu}\\
and\\
Carla D. Savage\\[-5pt]
\small Department of Computer Science, North Carolina State University,\\[-5pt]
\small Raleigh, NC 27695-8206, USA, {\tt savage@ncsu.edu}
}

\date{\today\\[10pt]
	\begin{flushleft}
	                                       \\[5pt]
	\end{flushleft}}

\maketitle

\begin{abstract}
We introduce the notion of a Mahonian pair.  Consider the set, $\bbP^*$, of all words having the positive integers as alphabet.  Given finite subsets $S,T\sbs\bbP^*$, we say that $(S,T)$ is a {\it Mahonian pair\/} if the distribution of the major index, $\maj$, over $S$ is the same as the distribution of the inversion number, $\inv$, over $T$.  So the well-known fact that $\maj$ and $\inv$ are equidistributed over the symmetric group, $\fS_n$, can be expressed by saying that $(\fS_n,\fS_n)$ is a Mahonian pair.  We investigate various Mahonian pairs $(S,T)$ with $S\neq T$. Our principal tool is Foata's fundamental bijection 
$\phi:\bbP^*\ra\bbP^*$ since it has the property that $\maj w = \inv \phi(w)$ for any word $w$.  We consider various families of words associated with Catalan and Fibonacci numbers.   We show that, when restricted to words in $\{1,2\}^*$, $\phi$ transforms familiar statistics on words into natural statistics on integer partitions such as the size of the Durfee square.
The Rogers-Ramanujan identities, the Catalan triangle, and various $q$-analogues also make an appearance.  We generalize the definition of Mahonian pairs to infinite sets and use this as a tool to connect a partition bijection of Corteel-Savage-Venkatraman  with the Greene-Kleitman decomposition of a Boolean algebra into symmetric chains.  We close with comments about future work and open problems.
\end{abstract}

\section{Introduction}

To introduce our principal object of study, Mahonian pairs, we need to set up some notation.
Let $\bbN$ and $\bbP$ denote the nonnegative and positive integers, respectively.  Consider the Kleene closure $\bbP^*$ of all words $w=a_1 a_2\ldots a_n$ where $a_i\in\bbP$ for $1\le i\le n$ and $n\ge0$. 
We let $l(w)=n$ be the length of $w$ and $\ep$ be the empty word of length $0$.  
We will often express words using multiplicity notation where $w=a^{m_1} b^{m_2} \ldots c^{m_k}$ is the word beginning with $m_1$ copies of $a$ followed by $m_2$ copies of $b$, and so forth.  Denote by $\Pi(w)$ the subset of $\bbP^*$ consisting of all permutations of $w$.

We will be considering various well-known statistics on $\bbP^*$.  The word $w=a_1 a_2\ldots a_n$ has {\it descent set\/}
$$
\Des w =\{i\ :\ \mbox{$1\le i<n$ and $a_i>a_{i+1}$}\},
$$
and {\it inversion set\/}
$$
\Inv w =\{(i,j)\ :\ \mbox{$1\le i < j\le n$ and $a_i>a_j$}\},
$$
and we say that $a_i$ is {\it in inversion with\/} $a_j$ if $(i,j)\in\Inv w$.
From these sets we get the {\it major index\/}
$$
\maj w =\sum_{i\in\Des w} i,
$$
the {\it descent number\/}
$$
\des w = |\Des w|,
$$
and the {\it inversion number\/}
$$
\inv w = |\Inv w|
$$
where $|\cdot|$ indicates cardinality.  For more information on these statistics, see Stanley's 
text~\cite[p.\ 21 \& ff.]{sta:ec1}.

Let $\fS_n\sbs\bbP^*$ be the symmetric group of all permutations of $\{1,\ldots,n\}$ whose elements will be viewed in one-line notation as sequences.  A celebrated result of MacMahon~\cite[pp.\ 508--549 and pp. 556--563]{mac:cp1} states that, for any $w\in\bbP^*$, the distribution of $\maj$ and $\inv$ over $\Pi(w)$ are the same, i.e., we have equality of the generating functions
\beq
\label{Mahonian}
\sum_{w\in\Pi(w)} q^{\maj w} = \sum_{w\in\Pi(w)} q^{\inv w}.
\eeq
A statistic $\sta:\Pi(w)\ra\bbN$ is called {\it Mahonian\/} if it is equidistributed with $\maj$ and $\inv$ and there is a large literature on such statistics both on words and other structures; see, for example, the work of Bj\"orner and Wachs~\cite{bw:qhl,bw:psl}.  This motivates our new definition.
Given finite subsets $S,T\sbe\bbP^*$, we call $(S,T)$ a {\it Mahonian pair\/} if 
\beq
\label{mp}
\sum_{w\in S} q^{\maj w} = \sum_{w\in T} q^{\inv w}
\eeq
So equation~\ree{Mahonian} can be expressed by saying that $(\Pi(w),\Pi(w))$ is a Mahonian pair.  

Note that the Mahonian pair  relation is not symmetric.  We wish to study various pairs where $S\neq T$.  Our main tool will be Foata's fundamental bijection~\cite{foa:nin} $\phi:\bbP^*\ra\bbP^*$ (defined in the next section) because it has the property that
\beq
\label{maj-inv-phi}
\maj w = \inv \phi(w)
\eeq
for all $w\in\bbP^*$.  So for any finite $S\sbe\bbP^*$ we have corresponding Mahonian pairs $(S,\phi(S))$ and $(\phi^{-1}(S),S)$.  Of course, the point of the definition~\ree{mp} is that $S$ and $T$ should have independent interest outside of being part of a pair.

The rest of this paper is structured as follows.  In the next section we recall the definition of Foata's bijection, $\phi$,  and prove some general results which will be useful in the sequel.  One of the goals of this work is to show that, when restricted to $\{1,2\}^*$, $\phi$ sends various well-known statistics on words to natural statistics on partitions; see, for example, equation~\ree{wSlaT}.
Section~\ref{bs} studies pairs associated with ballot sequences and Catalan numbers.    Pairs connected with Fibonacci numbers are studied in Section~\ref{wcf}.  The definition of Mahonian pairs is generalized to infinite sets in Section~\ref{imp} and we use this idea to  connect a bijection of Corteel-Savage-Venkatraman~\cite{csv:bpa} with the Greene-Kleitman decomposition of a Boolean algebra into symmetric chains~\cite{gk:svs}.  The final section is devoted to remarks and open problems.

\section{Foata's fundamental bijection}

We now review the construction and some properties of Foata's bijection $\phi:\bbP^*\ra\bbP^*$.  This section is expository and the results herein are not new unless otherwise stated.

Given $v=a_1 a_2\ldots a_n$ we inductively construct a sequence of words $w_1, w_2, \dots, w_n=\phi(v)$ as follows.  Let $w_1=a_1$.  To form $w_{i+1}$ from $w_i=b_1 b_2\ldots b_i$, compare $a_{i+1}$ with $b_i$.  Form the unique factorization $w_i = f_1 f_2 \cdots f_k$ such that, if $b_i\le a_{i+1}$ (respectively, $b_i>a_{i+1}$), then each factor contains only elements greater than (respectively, less than or equal to) $a_{i+1}$ except the last which is less than or equal to (respectively, greater than) $a_{i+1}$.  Let $g_j$ be the cyclic shift of $f_j$ which brings the last element of the factor to the front and let 
$w_{i+1}=g_1 g_2\cdots g_k a_{i+1}$.   As an example, finding $\phi(2121312)$ would give rise to the following computation where we separate the factors of $w_i$ with dots:
$$
\barr{lclcll}
w_1&=&2		&=&2					&\mbox{since $a_2=1$},\\
w_2&=&21		&=&2\cdot1				&\mbox{since $a_3=2$},\\
w_3&=&212		&=&2\cdot12				&\mbox{since $a_4=1$},\\
w_4&=&2211		&=&2\cdot2\cdot1\cdot1		&\mbox{since $a_5=3$},\\
w_5&=&22113	&=&2\cdot2\cdot113			&\mbox{since $a_6=1$},\\
w_6&=&223111	&=&2\cdot2\cdot31\cdot1\cdot1	&\mbox{since $a_7=2$},\\
w_7&=&2213112	&=&\phi(2121312).
\earr
$$
One can show that $\phi$ is a bijection by constructing a step-by-step inverse.  Equation~\ree{maj-inv-phi} can be shown to hold by noting that $w_i=\phi(v_i)$ where $v_i=a_1 a_2\ldots a_i$ and then showing that at each stage of the algorithm the change in $\maj$ in passing from $v_i$ to $v_{i+1}$ is the same as the change in $\inv$ in passing from $w_i$ to $w_{i+1}$.

We will be particularly interested in the action of $\phi$ on words $w\in\{1,2\}^*$.  Keeping the notation of the previous paragraph, note that if $a_{i+1}=2$ then all the factors are of length 1 and so $w_{i+1}=w_i 2$.  If $a_{i+1}=1$ then the factors will either be of the form $1^n2$ if $b_i=2$ or of the form $2^n1$ if $b_i=1$ for some $n\ge0$.  This enables us to give a nice characterization of $\phi$.  In the proof of this result and later on, it will be convenient to have another notation for words in $\{1,2\}^*$ by subscripting the ones right to left with $1,2,\ldots$ and the twos similarly left to right.  For example, $w=112122$  would become 
$w=1_3 1_2 2_1 1_1 2_2 2_3$.

The following recursive description of $\phi$ can also be obtained by combining Theorem~11.1 with equations~(11.8) and~(11.9) in the lecture notes of Foata and Han~\cite{fh:qsc}.

\ble
\label{philem}
The map $\phi$ on $\{1,2\}^*$ can be defined recursively by $\phi(\ep)=\ep$, $\phi(1)=1$ and the three rules
\ben
\item[{\rm (i )}] $\phi(w2) =\phi(w)  2$, 
\item[{\rm (ii )}] $\phi(w11) = 1 \phi(w1)$,
\item[{\rm (iii )}] $\phi(w21) = 2 \phi(w)  1$.
\een
\ele
\prf
Since the initial conditions and recursive rules uniquely define a map $\phi':\{1,2\}^*\ra\{1,2\}^*$, it suffices to show that $\phi$ satisfies these statements to show that $\phi'=\phi$.  The initial conditions and (i) follow directly from the description of $\phi$ given above.

(ii)  Let $v=\phi(w1)$.  To form $\phi(w11)$, one must cycle the factors of $v$ which are all of the form $2^n1$ for some $n$.  During this process, $1_i$ in $v$ moves into the position of $1_{i+1}$ for all $i$ except when $i$ has attained its maximum.  In that case, $1_i$ becomes a new one at the beginning of the resulting word.  And the final one in $w11$ takes the place of $1_1$ in $v$.  Thus the total effect is to move $v$ over one position and prepend a one, in other words $\phi(w11)=1\phi(w1)$.

(iii)  The proof is similar to that of (ii) while also using (i).  So it is left to the reader.
\hqedm

Note that $v\in\{1,2\}^*$ has $\des v = d$ if and only if one can write
\beq
\label{des}
v=1^{m_0} 2^{n_0} 1^{m_1} 2^{n_1} \ldots 1^{m_d} 2^{n_d}
\eeq
where $m_0, n_d\ge0$ and $m_i,n_j>0$ for all other $i,j$.  We will now derive a new, non-recursive description of $\phi$ on binary words which will be crucial to all that follows.
\bpr
\label{phi}
Let $v$ have $d$ descents and so be given by~\ree{des}.  It follows that
$$
\phi(v) =1^{m_d-1} 2 1^{m_{d-1}-1} 2 \ldots 1^{m_1-1} 2 1^{m_0} 2^{n_0-1} 1 2^{n_1-1} 1 \ldots 2^{n_{d-1}-1} 1 2^{n_d}.
$$
\epr
\prf
We induct on $d$.  When $d=0$ we have $v=1^{m_0} 2^{n_0}=\phi(v) $ which is correct.  For $d>1$, one can write
$$
v=u 2 1^{m_d} 2^{n_d}
$$
where $u$ is the appropriate prefix of $v$.  Using the previous lemma repeatedly gives
\bal
\phi(v)
   &=\phi(u 2 1^{m_d}) 2^{n_d}\\
   &=1^{m_d-1} \phi(u 2 1) 2^{n_d}\\
   &=1^{m_d-1} 2 \phi(u) 1 2^{n_d}.  
\end{align*}
By induction, this last expression coincides with the desired one.
\hqedm

Applying $\phi^{-1}$ to both sides of the last set of displayed equations, one gets the following result which we record for later use.
\ble
\label{phiinvlem}
If $w=1^m2u12^n$ for $m,n\ge0$ then

\vs{5pt}

$
\hfill
\phi^{-1}(w)=\phi^{-1}(u)21^{m+1}2^n.
\hfill\qed
$
\ele

\thicklines
\setlength{\unitlength}{2pt}
\bfi
$$
\bpi(60,50)(0,0)
\Gaaea \Gabeb \Gacec \Gaded \Gaeee \Gafef
\Gaaaf \Gbabf \Gcacf \Gdadf \Geaef
\multiput(0,0)(.1,0){11}{\line(0,1){20}}
\multiput(0,20)(0,-.1){11}{\line(1,0){21}}
\multiput(20,20)(.1,0){11}{\line(0,1){20}}
\multiput(20,40)(0,-.1){11}{\line(1,0){11}}
\multiput(30,40)(.1,0){11}{\line(0,1){10}}
\multiput(30,50)(0,-.1){11}{\line(1,0){10}}
\pmo{-3,3}{$1_5$}
\pmo{-3,13}{$1_4$}
\pmo{7,15}{$2_1$}
\pmo{17,15}{$2_2$}
\pmo{24,23}{$1_3$}
\pmo{24,33}{$1_2$}
\pmo{27,36}{$2_3$}
\pmo{34,44}{$1_1$}
\pmo{37,53}{$2_4$}
\epi
$$
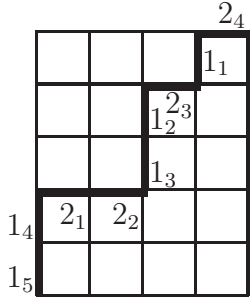
\capt{A partition $\la$ contained in $5\times 4$ and the corresponding lattice path}
\label{lafig}
\efi

There is a well-known intimate connection between words in $\{1,2\}^*$ and integer partitions.  An {\it integer partition\/} is a weakly decreasing sequence $\la=(\la_1,\la_2,\ldots,\la_k)$ of positive integers.  A cornucopia of information about partitions can be found in Andrews' book~\cite{and:tp}.   The $\la_i$ are called {\it parts\/} and we will use multiplicity notation for them  as we do in words.  Let $|\la|$  denote the sum of its parts. 
 The {\it Ferrers diagram of $\la$\/} consists of $k$ left-justified rows of boxes with $\la_i$ boxes in row $i$.  The Ferrers diagram for $\la=(3,2,2)=(3,2^2)$ is shown northwest of the dark path  in Figure~\ref{lafig}.   We let $(i,j)$ denote the box in row $i$ and column $j$.  We say that {\it $\la$  fits in an $m\times n$ rectangle\/}, written $\la\sbe m\times n$, if $\la\sbe(m^n)$ as Ferrers diagrams.  Figure~\ref{lafig} shows that $(3,2,2)\sbe 5\times 4$.
 The {\it Durfee square of $\la$\/}, $D(\la)$,  is the largest partition $(d^d)\sbe\la$.  We let $d(\la)=d$ denote the length of a side of $D(\la)$ and in our example $d(3,2,2)=2$.   In general, if $\mu\sbe\la$ then we have a {\it skew partition\/} $\la/\mu$ consisting of all the boxes of the Ferrers diagram of $\la$ which are not in $\mu$.  So the Ferrers diagram of $(4^5)/(3,2^2)$ consists of the boxes southeast of the dark path in the figure.  

Consider any word $w$ in the set $\Pi(1^m 2^n)$ of permutations of $1^m 2^n$.  We identify $w$ with a lattice path $P(w)$  in $\bbZ^2$ from $(0,0)$ to $(n,m)$, where a $1$ denotes a step one unit north and a $2$ a step one unit east.   Figure~\ref{lafig} displays the path for $w=112211212$ with each step labeled by the corresponding element of $w$ along with its subscript.
In this way we associate with $w$ an integer partition $\la(w)\sbe m\times n$ whose Ferrers diagram consists of the boxes inside the rectangle and northwest of $P(w)$.  Note  that $2_j$ is in inversion with $1_i$ if and only if the Ferrers diagram of $\la(w)$ contains the box $(i,j)$.   It follows that
\beq
\label{la}
\inv w = |\la(w)|.
\eeq  
Also, $d(\la)$ is the largest subscript such that $2_d$ is in inversion with $1_d$.

The following result will be useful in deriving generating functions using $\phi$.  That statement about $\maj$ is well known, but the equation involving $\des$ appears to be new.
\bco
\label{maj,des}
If $v\in\{1,2\}^*$ and $\la=\la(\phi(v))$, then
\ben
\item[{\rm (i)}] $\maj v= |\la|$, and
\item[{\rm (ii)}] $\des v =d(\la)$. 
\een
\eco
\prf
Let $w=\phi(v)$.  By the properties of Foata's map and~\ree{la} we have $\maj v=\inv w =|\la|$.

For (ii),  suppose $v$ has the form~\ree{des}.   From  Proposition~\ref{phi}  we can read off the positions of the first $d$ ones and the last $d$ twos in $w=\phi(v)$.  In particular, $2_d$ is in inversion with $1_d$ and this is not true for any larger subscript.  Thus 
$d(\la)=d=\des v$ as desired. 
\hqedm

\section{Ballot sequences}
\label{bs}

\subsection{Applying $\phi$}

We will now investigate the effects of $\phi$ and $\phi^{-1}$ on ballot sequences.
Say that $w\in\bbP^*$ is a {\it ballot sequence\/} if, for every prefix $v$ of $w$ and every $i\in\bbP$, the number of $i$'s in $v$ is at least as large as the number of $(i+1)$'s.
Consider 
$$
B_n=\{w\in\Pi(1^n 2^n)\ :\ \mbox{$w$ is a ballot sequence}\}.  
$$
The number of such ballot sequences is given by the {\it Catalan number\/}
\beq
\label{C_n}
C_n=\frac{1}{n+1}{2n \choose n}.
\eeq
There is a corresponding $q,t$-analogue
\beq
\label{c_n}
c_n(q,t)=\sum_{w\in B_n} q^{\maj w } t^{\des w }.
\eeq
The $c_n(q,t)$  were first defined in a paper of F\"urlinger and Hofbauer~\cite{fh:qcn}.  The case $t=1$ was the subject of an earlier note by Aissen~\cite{ais:vtp}.   These $q$-analogues have since been studied by various authors~\cite{ch:fgc,cww:ldc,kra:clp2}.  Note that one has the following $q$-analogue of~\ree{C_n},
\beq
\label{c_n(q,1)}
c_n(q,1)=\frac{1}{[n+1]}\gauss{2n}{n},
\eeq
where we have the usual conventions $[n]=1+q+\cdots+q^{n-1}$, $[n]!=[1][2]\cdots[n]$, and $\left[n\atop k\right]=[n]!/([k]![n-k]!)$.
No closed form expression is known for $c_n(q,t)$.

In order to describe $\phi(B_n)$, we need to define the ranks of a partition, $\la$.
We let the {\it $i$th rank\/} of $\la$, $1\le i\le d(\la)$,  be 
$$
r_i(\la) =\la_i-\la_i'
$$
where $\la'$ is the {\it conjugate of $\la$\/} obtained by transposing $\la$'s Ferrers diagram.    
We also let
$$
R_n =\{\la\ :\ \mbox{$\la\sbe n\times n$ and $r_i(\la)<0$ for all $1\le i\le d(\la)$ }\}.  
$$
The concept of rank goes back to Dyson~\cite{dys:sgt}.  Partitions with all ranks positive (which are in bijection with partitions with all ranks negative by conjugation) have arisen in the work of a number of authors~\cite{and:dsg,bs:rcg,er:gp,ra:ngp}.  Interest in them stems from connections with partitions which are the degree sequences of simple graphs and with the Rogers-Ramanujan identities. 
% It follows from the work of Corteel, Savage, and Venkatraman~\cite{csv:bpa} that $|R_n|= C_n$ which is also a corollary of the next result.
\bth
\label{phi(B_n)}
We have
\beq
\label{phi(B_n)eq}
\phi(B_n) = \{w\in\Pi(1^n 2^n)\ :\ \la(w)\in R_n\}.
\eeq
\eth
\prf
Assume that $v\in B_n$ has the form~\ree{des}.  Then $v\in B_n$ is equivalent to the fact that $m_0+\cdots+m_k\ge n_0+\cdots+ n_k$ for $0\le k\le d$ with both sides equalling $n$ for $k=d$.   

To check the ranks of $\la(w)$, note that from our description of the relationship between $w$ and $\la(w)$ it follows that $\la_i$ is the number of twos before $1_i$ in $w$, and $\la_j'$ is the number of ones after $2_j$.  Using Proposition~\ref{phi}, we have
\bal
\la_l&=d+(n_0-1)+(n_1-1)+\cdots+(n_{d-l}-1),\\
\la_l'&=d+m_0+(m_1-1)+\cdots+(m_{d-l}-1)
\end{align*}
for $1\le l\le d$.
Hence $r_l(\la(w))<0$ if and only if $\la_l'>\la_l$ for $1\le l\le d$, which is clearly equivalent to the ballot conditions for $0\le k<d$.  And we have the desired equality when $k=d$ since $w\in\Pi(1^n 2^n)$,
\hqedm

Combining the previous theorem with the definition of $c_n(q,t)$ and equation~\ree{c_n(q,1)}, we immediately obtain the following result.  The case $t=1$  was obtained by Andrews~\cite{and:dsg} using more sophisticated means.
\bco
We have
$$
\sum_{\la\in R_n} q^{|\la|} t^{d(\la)}=c_n(q,t).
$$
In particular,

\eqed{
\sum_{\la\in R_n} q^{|\la|} =\frac{1}{[n+1]}\gauss{2n}{n}.
}
\eco

\subsection{Applying $\phi^{-1}$}

We can obtain another Mahonian pair by applying $\phi^{-1}$ to $B_n$.  For the proof characterizing the preimage, we need a notion of conjugation for sequences.  If $w=b_1 b_2\ldots b_n\in\{1,2\}^*$ then let 
\beq
\label{prime}
w'=b_n'\ldots b_2' b_1'
\eeq 
where $b_i'=3-b_i$ for $1\le i\le n$, i.e., read the sequence backwards while exchanging the ones and twos.  It should be clear from the definitions that $\la(w)$ and $\la(w')$ are conjugate partitions.
\bth
\label{phi-1(B_n)}
Let $v$ have the form~\ree{des}.  Then $\phi(v) \in B_n$ if and only if we have the conditions:
\ben
\item[{\rm (i)}] for all $i$ with $1\le i\le d$
\bal
m_d+m_{d-1}+\cdots+m_{d-i+1}&\ge 2i,\\
n_d+n_{d-1}+\cdots+n_{d-i+1}&\ge 2i-1,
\end{align*}
\item[{\rm (ii)}] and
$$
\sum_{i=0}^d m_i =\sum_{i=0}^d n_i = n.
$$
\een
\eth
\prf
Suppose $\phi(v) \in B_n$.  Then, using Proposition~\ref{phi}, we see that the ballot sequence condition for the prefix up to $2_i$ (where $1\le i \le d$) is given by 
$$
(m_d-1) + (m_{d-1}-1) + \cdots + (m_{d-i+1} -1) \ge i.
$$
This, in turn, is clearly equivalent to the first inequality in (i).   

To obtain the second inequality, note that if $w$  is in  $B_n$ then so is $w'$ as defined in equation~\ree{prime}.  Now use the same reasoning as in the previous paragraph.  

Finally, condition (ii) follows since $\phi$ preserves the number of ones and twos.  It is easy to see that the reasoning above is reversible and so (i) and (ii) imply $\phi(v) \in B_n$.
\hqedm

Let us put the previous theorem in context.  Consider
$$
B_{m,n}=\{w\in\Pi(1^m 2^n)\ :\ \mbox{$w$ is a ballot sequence}\}
$$
The {\it Catalan triangle\/} (Sloane A008315) has entries $C_{n,d}$ for $0\le d \le \fl{n/2}$, where
$$
C_{n,d} = |B_{n-d,d}|.
$$
It is well known that
\beq
\label{sq}
C_n=\sum_{d\ge0} C_{n,d}^2
\eeq
where we assume $C_{n,d}=0$ for $d>\fl{n/2}$.  One way to see this is to consider the map 
$\be:B_n\ra \uplus_{d\ge0} B_{n-d,d}^2$ defined as follows.  Given $w=b_1 b_2\ldots b_{2n}\in B_n$, write $w=xy$ where $x=b_1 b_2\ldots b_n$ and $y=b_{n+1} b_{n+2}\ldots b_{2n}$, and define
\beq
\label{be}
\be(xy)=(x,y')
\eeq
where $y'$ is as in~\ree{prime}.
It is easy to show that $\be$ is well defined and a bijection and so~\ree{sq} follows.

Associated with any $v$ of the form~\ree{des} are two compositions (ordered partitions allowing zeros), the {\it one's composition\/} $\om(v) =(m_0,m_1,\ldots,m_d)$ and the {\it two's composition\/} $\tau(v) =(n_0,n_1,\ldots,n_d)$.  Note that $\om(v) \in \bbN \bbP^d$ while 
$\tau(v) \in \bbP^d\bbN$.  Consider two sets of compositions
$$
O_{n,d}=\{\om\in\bbN\bbP^d\ :\ \makebox{$l(\om)= n$ and $\om_d+\om_{d-1}+\cdots+\om_{d-i+1}\ge 2i$ for $1\le i\le d$}\},
$$
and
$$
T_{n,d}=\{\tau\in\bbP^d\bbN\ :\ \makebox{$l(\tau)= n$ and $\tau_d+\tau_{d-1}+\cdots+\tau_{d-i+1}\ge 2i-1$ for $1\le i\le d$}\}.
$$

To describe the next result, let $A_{n,d}$ be the set of $v$ of the form~\ree{des} such that $\phi(v) \in B_n$.  Also note that if we define a map $f$ on $A_{n,d}$ for each $d\ge0$, then $f$ can  be considered as a map on 
$\uplus_{d\ge0} A_{n,d}=\phi^{-1}(B_n)$.  Similar considerations apply to the other sets with subscripts $n,d$ defined above.
\bth
We have the following facts.
\ben
\item[{\rm(i)}] The maps $o:O_{n,d}\ra B_{n-d,d}$ and $t:T_{n,d}\ra B_{n-d,d}$ defined by
$$
o(\om_0,\om_1,\ldots,\om_d)=1^{\om_d-1} 2 1^{\om_{d-1}-1} 2 \ldots 1^{\om_1-1} 2 1^{\om_0},
$$ 
and
$$
t(\tau_0,\tau_1,\ldots,\tau_d)= 1^{\tau_d} 2 1^{\tau_{d-1}-1} 2 1^{\tau_{d-2}-1}  \ldots 2 1^{\tau_0-1} 
$$
are bijections.  Thus
$$
|O_{n,d}| = |T_{n,d}| = C_{n,d}.
$$
\item[{\rm(ii)}] The map $\om\times\tau:A_{n,d}\ra O_{n,d}\times T_{n,d}$ given by $v\mapsto (\om(v),\tau(v))$ is a bijection.  Thus
$$
|A_{n,d}|= C_{n,d}^2.
$$
\item[{\rm(iii)}]  Composing maps from right to left, we have
$$
\be=(o\times t) \circ (\om\times\tau) \circ \phi^{-1}.
$$
where $\be$ is the map given by~\ree{be}.
\een
\eth
\prf
(i)  We will prove the statements involving $O_{n,d}$ as proofs for $T_{n,d}$ is similar.  
We must first show that $o$ is well defined, i.e., that if 
$\om\in O_{n,d}$ then $o(\om)\in B_{n-d,d}$.  The number of ones in $o(\om)$ is $\om_0+\sum_{i=1}^d (w_i-1) = n-d$ as desired, and it is clear that there are exactly $d$ twos.  To verify the ballot condition, it suffices to check the prefix ending in $2_j$ for $1\le j\le d$.  Using the defining inequality for $O_{n,d}$, the number of ones in this prefix is
$$
\sum_{i=1}^j (\om_{d-i+1}-1) \ge 2j - j = j
$$
which is what we need.  Constructing an inverse to prove bijectivity is easy, and the statement about cardinalities follows from $o$ being bijective.

\medskip

(ii) The fact that the map is well defined and bijective is just a restatement of Theorem~\ref{phi-1(B_n)}.  The cardinality of $A_{n,d}$ can now be computed using this bijection and (i).

\medskip

(iii)  Any $w\in B_n$ can be written uniquely as $w=xy$ where $l(x)=l(y)=n$.  Say $x\in\Pi(1^{n-d}2^d)$ for some $d\ge0$.  It follows that $y\in\Pi(1^d 2^{n-d})$. So we can write
\bal
x&= 1^{k_d} 2 1^{k_{d-1}} 2\ldots 1^{k_1} 2 1^{k_0},\\
y&= 2^{l_0} 1 2^{l_1} 1\ldots 2^{l_{d-1}} 1 2^{l_d}
\end{align*}
where $k_0,l_0\in\bbN$ and all the other $k_i$ and $l_j$ are positive.  Thus we have $\be(w) =(x,y')$ where 
$$
y'= 1^{l_d} 2 1^{l_{d-1}} 2 \ldots 1^{l_1} 2 1^{l_0}.
$$

On the other side of the desired equality, we first use Lemma~\ref{phiinvlem} to compute
\bal
w&\stackrel{\phi^{-1}}{\mapsto} 1^{k_0} 2^{l_0+1} 1^{k_1+1} 2^{l_1+1} \ldots 2^{l_{d-1}+1} 1^{k_d+1} 2^{l_d}\\
   &\stackrel{\om\times\tau}{\mapsto} \left((k_0,k_1+1,k_2+1,\ldots,k_d+1),\ (l_0+1,l_1+1,\ldots,l_{d-1}+1,l_d)\right)\\
   &\stackrel{o\times t}{\mapsto} (1^{k_d} 2 1^{k_{d-1}} 2\ldots 1^{k_1} 2 1^{k_0},\ 
									 1^{l_d} 2 1^{l_{d-1}} 2 \ldots 1^{l_1} 2 1^{l_0}).
\end{align*}
This agrees exactly with the image of $\be$ obtained in the previous paragraph.
\hqedm

Note that~\ree{sq} follows as a corollary to part (ii) of this theorem.  It is also possible to give  $q$-analogues of this equality.  Let
$$
C_{n,d}(q)=\sum_{w\in B_{n-d,d}} q^{\inv w}
$$
with $C_n(q)=C_{2n,n}(q)$.   These $q$-Catalan numbers were first defined by Carlitz and Riordan~\cite{cr:tel}.  They have since become the objects of intense study. This is in part because of their connection with the space of diagonal harmonics;  see the text of Haglund~\cite{hag:qtc} for more information and references.

We also let
$$
c_{n,d}(q,t)=\sum_{w\in B_{n-d,d}} q^{\maj w } t^{\des w }.
$$
Comparing with~\ree{c_n} gives $c_n(q,t)=c_{2n,n}(q,t)$.  We will often use upper and lower case letters for the generating functions of $\inv$ and $\maj$, respectively, over the same set.
It will also be convenient to have the notation
$$
\de_{n,d}(q,t)=c_{n,d}(q,t)-c_{n-1,d-1}(q,t).
$$

Before stating the $q$-analogues, it will be useful to have a result about how our statistics change under the prime operation.
\ble
\label{y'}
Suppose $l(y)=n$.
\ben
\item[{\rm(i)}] $\inv y'=\inv y  $.
\item[{\rm(ii)}] $\des y'   =\des y  $.
\item[{\rm(iii)}] $\maj y'   = n\des y  -\maj y  $.
\een
\ele
\prf
Let $y=a_1 a_2\ldots a_n$.  Note that $(i,j)$ is an inversion of $y$ if and only if $i<j$ and $a_i'<a_j'$, which is if and only if $(n-j+1,n-i+1)$ is an inversion of $y'$.  All three parts of the lemma now follow.
\hqedm

\bpr
We have
$$
C_n(q)=\sum_{d\ge0} q^{d^2} C_{n,d}(q)^2,
$$
and
$$
\barr{l}
c_n(q,t)=\dil\sum_{d\ge0}\left[
q^n t c_{n-1,d-1}(q,t) c_{n-1,d-1}(q^{-1},q^{2n}t)+c_{n-1,d-1}(q,t) \de_{n,d}(q^{-1},q^{2n}t)\right.\\[20pt]
\left.\hs{65pt} +\de_{n,d}(q,t) c_{n-1,d-1}(q^{-1},q^{2n}t)+\de_{n,d}(q,t) \de_{n,d}(q^{-1},q^{2n}t)
\right].
\\[10pt]
\earr
$$
\epr
\prf
Suppose $w=b_1b_2\ldots b_{2n}\in B_n$ and let $x=b_1\ldots b_n$, $y=b_{n+1}\ldots b_{2n}$ so that  $\be(w) =(x,y')$. 
If $x\in B_{n-d,d}$ then $y\in B_{d,n-d}$ and, by Lemma~\ref{y'}(i),
$$
\inv w  = \inv x +\inv y   + d^2 = \inv x +\inv y'   +d^2
$$
where the $d^2$ comes from the inversions between twos in $x$ and ones in $y$.   Using the fact that $\be$ is a bijection, the displayed equation above translates directly into the first desired formula.

To obtain the second formula, it will be convenient to have the characteristic function $\chi$ which is $1$ on true statements and $0$ on false ones.  Keeping the notation from the first paragraph, we have
\bal
\des w  &= \des x  + \des y   + \chi(b_n>b_{n+1}),\\
\maj w  &= \maj x  + \maj y   + n\des y   + n\chi(b_n>b_{n+1})
\end{align*}
where the $n\des y  $ term comes from the fact that $i$ is a descent of $y$ if and only if $n+i$ is a descent of $w=xy$.  Note also that $b_n>b_{n+1}$ is equivalent to $b_n=2$ and $b_{n+1}=1$.  Using Lemma~\ref{y'}(ii) and (iii), we obtain
\bal
\des w  &= \des x  + \des y'    + \chi(b_n=b_{n+1}'=2),\\
\maj w  &= \maj x  - \maj y'    + 2n\des y'    + n\chi(b_n=b_{n+1}'=2).
\end{align*}

The proof now breaks down into four cases depending on the values of $b_n$ and $b_{n+1}$.  These correspond to the four terms in the expression for $c_n(q,t)$.  Since  the proofs are similar, we will only derive the second of them which corresponds to the case $b_n=2$ and $b_{n+1}'=1$.  Using the last pair of displayed equations, we obtain in this case
$$
\sum_w q^{\maj w }t^{\des w }
=  \sum_x q^{\maj x } t^{\des x } \sum_{y'} q^{-\maj y'   } (q^{2n} t)^{\des y'   }.
$$
Now $x$ is ranging over all elements of $B(n-d,d)$ which end in a two.  But this two does not contribute to $\maj$ or $\des$.  So removing the two, we get the same sum over all elements of $B(n-d,d-1)$.  Thus this sum can be replaced by $c_{n-1,d-1}(q,t)$.  The $y'$ sum  ranges over elements of $B(n-d,d)$ which end in a $1$.  Rewriting this as the sum over all of $B(n-d,d)$ minus the sum over those elements ending in a $2$, and then using  arguments similar to those concerning the $x$ sum,  gives us the $ \de_{n,d}(q^{-1},q^{2n}t)$ factor.
\hqedm

%CITE AND COMPARE THE FIRST $q$-analogue with $q$-Vandermonde as in Andrews' book p.\ 37 (33.10) where $m=n=h$ and %$k\ra n-k$.

Several observations about the expression for $c_n(q,t)$ are in order.
Note that, using the definition of $\de_{n,d}(q,t)$,  one can simplify the sum by combining either the 2nd and 4th or the 3rd and 4th terms.  But we have chosen the more symmetric form.   Note also that such an expression can not be reduced to one in $q$ alone because of the appearances of $q^{2n}$ in the second variable.  

\subsection{Extensions}

There are analogues of Theorems~\ref{phi(B_n)} and~\ref{phi-1(B_n)} for ballot sequences with any number of ones and twos.  Let
$$
R_{k,l} =\{\la\ :\ \mbox{$r_i(\la)<0$ for $1\le i\le d(\la)$ and $\la\sbe k\times l$}\}.  
$$
We will just state the next results as their proofs are similar to the case when $k=l$.
\bth
Suppose $k\geq l$.  We have
$$
\phi(B_{k,l})=\{w\in\Pi(1^k2^l)\ :\ \la(w)\in R_{k,l}\}.
$$
\eth

\bth
For $k\geq l$ and 
 for $v$ having the form~\ree{des}, we have $\phi(v) \in B_{k,l}$ if and only if 
\ben
\item[{\rm(}i{\rm )}] for all $i$ with $1\le i\le d$
\bal
m_d+m_{d-1}+\cdots+m_{d-i+1}&\ge 2i,\\
n_d+n_{d-1}+\cdots+n_{d-i+1}+(k-l)&\ge 2i-1,
\end{align*}
\item[{\rm (}ii{\rm )}] as well as

$
\hfill\dil\sum_{i=0}^d m_i =k, \qmq{and}
\dil\sum_{i=0}^d n_i = l.\hfill\qed
$
\een
\eth

There is another way to generalize Theorem~\ref{phi(B_n)}.  Given $w\in\{1,2\}^*$ we let $e(w)$ be the maximum excess of the number of twos over the number of ones in any prefix of $w$.  So, if $w$ has the form~\ree{des}, then
let
$$
e_i(w)=n_0+n_1+\cdots + n_i - m_0 - m_1 -\cdots - m_i
$$
and we have
$$
e(w)=\max_{0\le i\le d} e_i.
$$
Note that $w$ is a ballot sequence if and only if $e(w)\le0$.  There is another combinatorial interpretation of $e(w)$ as follows.  One can pair up certain ones and twos (thinking of them as left and right parentheses, respectively) of $w$ in the usual manner:  If a one is immediately followed by a two then they are considered paired.  Remove all such pairs from $w$ and iteratively pair elements of the remaining word.  Let $p(w)$ be the number of pairs in $w$.  If $w\in\Pi(1^n2^n)$, then it is easy to see by induction on $n$ that
$$
e(w) = n - p(w).
$$
We will return to pairings in Section~\ref{imp}.
The statistic on integer partitions $\la$ corresponding to $e(w)$ is the maximum rank
$$
r(\la)=\max_{1\le i\le d(\la)} r_i(\la).
$$
\ble
\label{e>=r+1}
Suppose $v\in\{1,2\}^*$, $\la=\la(\phi(v))$, and $d=d(\la)$.  It follows that
$$
e(v) \ge r(\la) +1.
$$
If $e_d(v)< e(v)$ (in particular, if $v\in\Pi(1^n2^n)$ is not a ballot sequence) then we have equality.
\ele
\prf
As usual, let $v$ have the form~\ree{des}.  
The same reasoning as in the proof of Theorem~\ref{phi(B_n)} shows that
\beq
\label{e_i(v)}
e_i(v)= r_{d-i}(\la)+1
\eeq
where $d=\des v = d(\la)$  and $0\le i<d$.  It follows that
$$
r(\la)+1 = \max_{0\le i<d} e_i(v) \le e(v)
$$
giving the desired inequality.  It should also be clear why $e_d(v)< e(v)$ implies equality.  Finally, suppose  $v\in\Pi(1^n2^n)$ is not a ballot sequence.   It follows that for some $i>d$ we have $e_i(v)>0=e_d(v)$ and we are done.
\hqedm

To state the analogue of Theorem~\ref{phi(B_n)}, let
$$
E_{n,k}=\{w\in\Pi(1^n2^n)\ :\ e(w)=k\}
$$
and
$$
P_{n,k}=\{w\in\Pi(1^n 2^n)\ :\ r(\la(w))=k\}.
$$
As an immediate corollary of the previous lemma we have the following.
\bth
If $k>0$, then the pair $(E_{n,k},P_{n,k-1})$ is Mahonian.\hqed
\eth

\section{Words counted by Fibonacci numbers}
\label{wcf}

We will explore Mahonian pairs constructed from various sets of sequences enumerated by the Fibonacci numbers.  Let the Fibonacci numbers themselves be defined by $F_0=F_1=1$ and $F_n=F_{n-1}+F_{n-2}$ for $n\ge2$.

The first set of sequences is given by
$$
\cF_n=\{w\in\{1,2\}^*\ :\ \mbox{$l(w)=n$ and $w$ has no consecutive ones}\}.
$$ 
It is well known and easy to show that $|\cF_n|=F_{n+1}$.   Let
$$
f_n(q,t)=\sum_{w\in\cF_n} q^{\maj w} t^{\des w}.
$$
As with the $F_n$, we have a recursion 
$$
f_n(q,t)=f_{n-1}(q,t)+q^{n-1} t  f_{n-2}(q,t)
$$
since any $w\in \cF_n$ can be obtained by appending $2$ to an element of $\cF_{n-1}$ (which does not change $\maj$ or $\des$) or by appending  $21$ to an element of $\cF_{n-2}$ (which increases $\maj$ by $n-1$ and $\des$ by $1$).  Polynomials satisfying this same recursion were introduced by Carlitz~\cite{car:qfn,car:qfp} and since studied by a number of authors, see the paper of Goyt and Sagan~\cite{gs:sps} for a comprehensive list.  However, the initial conditions used for these polynomials is different from ours, and so a different sequence of polynomials is generated.  No literature seems to exist about the other $q$-analogues of $F_n$ given by taking the distribution of $\inv$ or $\maj$ over the various sets considered in this section.   It would be interesting to study their properties.

To state our result about what $\phi$ does to $\cF_n$, it will  be convenient to consider
$$
\cF_{n,k}=\{w\in\cF_n\ :\ \mbox{$w$ has $k$ ones}\}.
$$
We also note that if partition $\la$ has Durfee square $D(\la)$ then $\la=D(\la)\uplus R(\la) \uplus B(\la)$ where $R(\la)$ and $B(\la)$ are the connected components of the skew partition $\la/D(\la)$ to the right and below $D(\la)$, respectively.  Finally, we let 
$\left[n\atop k \right]=0$ if $k<0$ or $k>n$.
\bth
\label{phi(cF)}
We have 
$$
\phi(\cF_{n,k})=\{ w\ :\  \mbox{$\la(w)=(\la_1,\ldots,\la_k)$ with $\la_1\le n-k$ and $\la_k\ge k-1$}\}.
$$
It follows that
\beq
\label{f_n}
f_n(q,t)=\sum_{k\ge0} q^{k(k-1)} t^{k-1} \left(\gauss{n-k}{k-1}+q^k t\gauss{n-k}{k}\right).
\eeq
\eth
\prf
Suppose $v\in\cF_{n,k}$ and let $\la(\phi(v))=\la$.  Since $v$ has $k$ ones and $n-k$ twos, $\la$ will fit in a $k\times(n-k)$ rectangle and so $\la_1\le n-k$.  Suppose $v$ has the form~\ree{des}.  So $v$ having no consecutive ones implies $m_i=1$ for $1\le i\le d=\des v$.  It follow from Proposition~\ref{phi} that 
$$
\la_k = \mbox{(the number of initial twos in $\phi(v)$)}\ge d \ge k-1
$$
where the last inequality follows from the fact that, since the ones in $v$ are not consecutive,  each of them except possibly the last creates a descent.  This reasoning is reversible and thus we have our characterization of $\phi(\cF_{n,k})$.

To get~\ree{f_n}, first use Corollary~\ref{maj,des} to write
$f_n(q,t)=\sum_\la q^{|\la|} t^{d(\la)}$ where the sum is over all $\la$ in the description of $\phi(\cF_{n,k})$ for each $k$.  Since $\la$ has $k$ parts and $\la_k\ge k-1$, $d(\la)=k-1$ or $d(\la)=k$, corresponding to the two terms in the summation.  If $d(\la)=k-1$ then we will have a factor of $q^{(k-1)^2}t^{k-1}$.  Also in this case $B(\la)=(k-1)$ and $R(\la)\sbe (k-1)\times (n-2k+1)$ because $\la\sbe k\times(n-k)$.  This gives a second factor of $q^{k-1}\left[n-k \atop k-1\right]$, and multiplying the two factors gives the first term in~\ree{f_n}.  The second term is obtained similarly when $d(\la)=k$.
\hqedm

Note that letting $t=1$ in~\ree{f_n} and applying one of the usual recursions for the $q$-binomial coefficients yields
$$
f_n(q,1)=\sum_{k\ge0} q^{k(k-1)} \gauss{n-k+1}{k}.
$$
Now letting $n\ra\infty$ gives, on the right hand side, $\sum_{k\ge0} q^{k^2-k}/(q)_k$  where we are using the Pochhammer symbol $(q)_k=(1-q)(1-q^2)\cdots(1-q^k)$.  This series was first studied by Carlitz~\cite{car:nrr,car:frr} who related it to the Rogers-Ramanujan identities.  Garrett, Ismail, and Stanton~\cite{gis:vrr} generalized the Rogers-Ramanujan identities to sums of the form
$\sum_{k\ge0} q^{k^2+mk}/(q)_k$ for any nonnegative integer $m$  from which one can also easily derive formulas for negative $m$ (so $m=0,1$ are the original identities and $m=-1$ is the case Carlitz considered).

To describe $\phi^{-1}(\cF_n)$ we will need a few definitions.  A {\it run\/}  in a word $w=b_1b_2\ldots b_n$ is a maximal factor $r=b_ib_{i+1}\ldots b_j$ such that $b_i=b_{i+1}=\ldots=b_j$.  We call $r$ a {\it $k$-run\/} if the common value of its elements is $k$.  Furthermore, $r$ is called the {\it prefix\/} or {\it suffix\/} run if $i=1$ or $j=n$, respectively.  For example, $w=1112212222$ has four runs, the prefix run $111$, the suffix run $2222$, the $1$-run consisting of the rightmost $1$, and the $2$-run $22$.
\bth
\label{phi-1(cF)}
The set $\phi^{-1}(\cF_n)$ consists of all $v\in\{1,2\}^*$ of length $n$ satisfying the following two conditions.
\ben
\item[{\rm(i)}]  If $v$ has a $1$-run $r$ as a prefix or suffix then $l(r)\le 1$ or $l(r)\le 2$, respectively.
\item[{\rm(ii)}] For any run of $v$ which is neither the prefix nor the suffix, we have $l(r)\le 2$ for $1$-runs and $l(r)\ge2$ for $2$-runs.
\een
\eth
\prf
Assume $v\in\phi^{-1}(\cF_n)$ has form~\ree{des}.  By Proposition~\ref{phi}, $\phi(v)$ has no consecutive ones if and only if the following three conditions hold
\ben
\item $m_i-1\le 1$ for $1\le i\le d$, 
\item $n_j-1\ge 1$ for $1\le j <d$, and 
\item $m_0\le 1$ with $n_0-1\ge 1$ if $m_0=1$. 
\een
It is easy to see that these are equivalent to the run conditions.
\hqedm

Dually to $\cF_n$, one can consider
$$
\cG_n=\{w\in\{1,2\}^*\ :\ \mbox{$l(w)=n$ and $w$ has no consecutive twos}\}.
$$
It is a simple matter to adapt the techniques use to demonstrate Theorems~\ref{phi(cF)} and~\ref{phi-1(cF)} to prove analogous results about $\cG_n$.  So we leave the details to the reader.

A third set of sequences counted by $F_n$ is
$$
\cH_n=\{w=b_1b_2\ldots\in\{1,2\}^*\ :\ \sum_i b_i =n\}.
$$
It is easy to prove that $|\cH_n|=F_n$ and this interpretation of the Fibonacci numbers can be made geometric by using tilings.  See the book of  Benjamin and Quinn~\cite[p.\ 1]{bq:prc} for examples.  Since $\phi$ preserves the number of ones and twos, it is clear that $\phi(\cH_n)=\cH_n$.  We have proved the following result.
\bth
The pair $(\cH_n,\cH_n)$ is Mahonian.  In other words

\vs{5pt}

$
\hfill\dil \sum_{w\in\cH_n} q^{\maj w}=\sum_{w\in\cH_n} q^{\inv w}.\hfill\qed
$
\eth

\section{Infinite Mahonian pairs}
\label{imp}

It is an easy matter to generalize the definition of a Mahonian pair to sets of any cardinality.  Call pair $(S,T)$ with $S,T\sbe\bbP^*$ {\it Mahonian\/} if there is a bijection $\al:S\ra T$ such that $\maj v = \inv \al(v)$ for all $v\in S$.  In this section we will study infinite Mahonian pairs connected to ballot sequences.  In particular, we will be able to show that a bijection of Corteel-Savage-Venkatraman~\cite{csv:bpa} between partitions with all ranks positive and partitions with no part equal to $1$ is essentially conjugation by $\phi$ of the map used by  Greene-Kleitman to obtain a symmetric chain decomposition of the Boolean algebra~\cite{gk:svs}.

Our starting point is an observation made independently by Andrews~\cite{and:dsg}  and Erd\"os-Richmond~\cite{er:gp}.  Let
$P$ be the set of all integer partitions.  Define
$$
R_{\ge t}=\{\la\in P\ :\  \mbox{$r_i(\la)\ge t$ for all $i$}\}
$$
and similarly for $R_{\le t}$.  Also, let
$$
P_{\neq t}=\{\la\in P\ :\ \mbox{$\la_i\neq t$ for all $i$}\}.
$$
\bth[\cite{and:dsg,er:gp}] 
\label{and-er}
We have

\vs{5pt}

$
\hfill 
\dil\sum_{\la\in R_{\ge1}} q^{|\la|} = \sum_{\la\in P_{\neq1}} q^{|\la|}.
\hfill\qed
$
\eth

We wish to make a couple of comments about this theorem.  First of all, it is well known that
\beq
\label{P_neq1}
 \sum_{\la\in P_{\neq1}} q^{|\la|}=\prod_{i\ge2} \frac{1}{1-q^i}.
\eeq
Secondly,  Theorem~\ref{and-er} is a special case of an earlier result of Andrews~\cite{and:ste} as generalized by Bressoud~\cite{bre:eps}.  (Andrews had the restriction that $M$ below must be odd.)
\bth[\cite{and:ste,bre:eps}]
\label{and-bre}
Let $M$, $r$ be integers satisfying $0<r<M/2$.
The number of partitions of $n$  whose ranks lie in the
interval $[-r+2, M-r-2]$ equals 
the number of partitions of $n$ with no part congruent to 0 or $\pm r$
modulo $M$ .\hqed
\eth

Note that if $M=n+2$ and $r=1$ then one recovers Theorem~\ref{and-er}.  And for $M=5$ and $r=1,2$ one obtains the Rogers-Ramanujan identities.

In order to involve Foata's map, we will need to associate with each partition $\la$  the word $w(\la)$ gotten by recording the north and east steps along the southeast boundary of $\la$'s Ferrers diagram with ones and twos, respectively.  To illustrate, $w(3,2,2)=221121$ as can be seen from Figure~\ref{lafig}.  Clearly $w=w(\la)$ for some nonempty $\la$ if and only if $w\in 2\{1,2\}^*1$, i.e., $w$ begins with a two and ends with a one.  We will often abuse notation and write things like $\phi(\la)$ for the more cumbersome $\phi(w(\la))$.  

To state our first result about infinite Mahonian pairs, we will use the notation
$$
W_v = \{1,2\}^* v\uplus\{\ep\}
$$
for any word $v$, and
$$
B_v=\{w\in W_v\ :\ \mbox{$w$ is a ballot sequence}\}.
$$
Furthermore, for any set $\La$ of partitions we let $\La'=\{\la'\ :\ \la\in \La\}$.  In the generating functions of equation~\ref{wSlaT} below, the exponent of $z$ may be negative.  So this should be viewed as an equality of Laurent series.  Finally, since the sets involved are infinite, one must be careful that the sums converge as formal Laurent series.  This can be seen by considering the right-hand side since there are only finitely many partitions with a given $|\la|$.
\bth
We have
$$
\phi^{-1}(P)=W_{21},\quad
\phi^{-1}(R_{\ge1}')=B_{21},\qmq{and}\phi^{-1}(P_{\neq1}')=W_{121}.
$$
It follows that $(S,T)=(W_{21},P)$, $(B_{21},R_{\ge1}')$, and $(W_{121},P_{\neq1}')$ are Mahonian pairs.  Furthermore, for any of these pairs $(S,T)$ we have
\beq
\label{wSlaT}
\sum_{v\in S} q^{\maj v} t^{\des v} z^{e(v)} = \sum_{\la\in T} q^{|\la|} t^{d(\la)} z^{r(\la)+1}.
\eeq
We also have
$$
\sum_{w\in B_{21}} q^{\maj w} =\sum_{w\in W_{121}} q^{\maj w} = \prod_{i\ge2}\frac{1}{1-q^i}.
$$
\eth
\prf
The fact that $\phi^{-1}(P)=W_{21}$ follows from Lemma~\ref{phiinvlem} and the observation above that $w=w(\la)$ precisely when $w\in 2\{1,2\}^*1$.   Note that if $w\in W_{21}$ then $e_d(w)<e_{d-1}(w)$, and this implies equality in Lemma~\ref{e>=r+1}.  Combining this with Corollary~\ref{maj,des} yields~\ree{wSlaT} for the case $(W_{21},P)$.  Since the other two pairs are formed from subsets of $W_{21}$ and $P$, the same equation will hold for them once we show that they are mapped in the desired way by $\phi$.

To prove that $\phi^{-1}(R_{\ge1}')=B_{21}$, note that $R_{\ge1}'=R_{\le-1}$.  Let us demonstrate that $\phi(B_{21})\sbe R_{\le-1}$.  Suppose that $v\in B_{21}-\{\ep\}$ is of the form~\ree{des} and let $w=\phi(v)$.   Then by Lemma~\ref{philem}(iii) we have $\phi(v)\in 2\{1,2\}^* 1$ and so there is a partition $\la$ with $w(\la)=w$.  Since $v$ is a ballot sequence we have $e(v)\le0$.  It follows from Lemma~\ref{e>=r+1} that $r(\la)\le -1$ giving the desired set containment.

For the reverse containment, start with a nonempty partition $\la\in R_{\le-1}$ and let $w=w(\la)$.  From the first paragraph of the proof we already know that  $v=\phi^{-1}(\la)\in W_{21}$ and this implies equality in Lemma~\ref{e>=r+1}. So, since $\la\in R_{\le-1}$, we have
$e(v)=r(\la) + 1 \leq 0$,
which is equivalent to $v$ being a ballot sequence.  Thus $v\in B_{21}$ which is what we needed to prove.

The demonstration of $\phi^{-1}(P_{\neq1}')=W_{121}$ is similar to the one just given, using the fact that $\la\in P_{\neq1}'$ if and only if $w=w(\la)\in 2\{1,2\}^*11$.  The final statement about generating functions now follows from  Lemma~\ref{y'}(i),  Theorem~\ref{and-er}, and equation~\ree{P_neq1}.
\hqedm

Corteel, Savage, and Venkatraman~\cite{csv:bpa} gave a bijective proof of Theorem~\ref{and-er}.  Using $\phi$, we can relate their function to the Greene-Kleitman symmetric chain decomposition of a Boolean algebra~\cite{gk:svs}.  For simplicity in comparing with our results thus far, we will describe the conjugate of their map.  For a partition $\la$ it will be convenient to define $\de(\la)=\la_1-\la_2$ and let
$$
D_t=\{\la\in P\ :\ \de(\la)=t\}.
$$
Note that $P_{\neq1}'=D_0$.  Define a map $CSV:D_0\ra R_{\le-1}$ by the following algorithm:

\ben
\item[CSV1]  Input $\la\in D_0$.
\item[CSV2] While $r(\la)\ge 0$ do
\ben
\item  Let $i$ be the maximum index such that $r_i(\la)=r(\la)$.
\item  Remove a part of size $i$ from $\la'$, add a part of size $i-1$ to $\la$, and increase the size of $\la_1$ by one.
\een
\item[CSV3] Output $\la$.
\een

We have written out the complete algorithm for computing $CSV(8,8,6,5,2,1)$ below  where dots rather than boxes have been used for the Ferrers diagrams. 
At each stage, the rank vector $\rho(\la)=[r_1(\la),\ldots,r_d(\la)]$ is displayed along with the maximum $r=r(\la)$ and the largest index $i=i(\la)$ where $r$ is achieved.  For future reference, we have also displayed $w=w(\la)$, $v=\phi^{-1}(w)$, and $\vep(v)=[e_0(v),\ldots,e_d(v)]$.

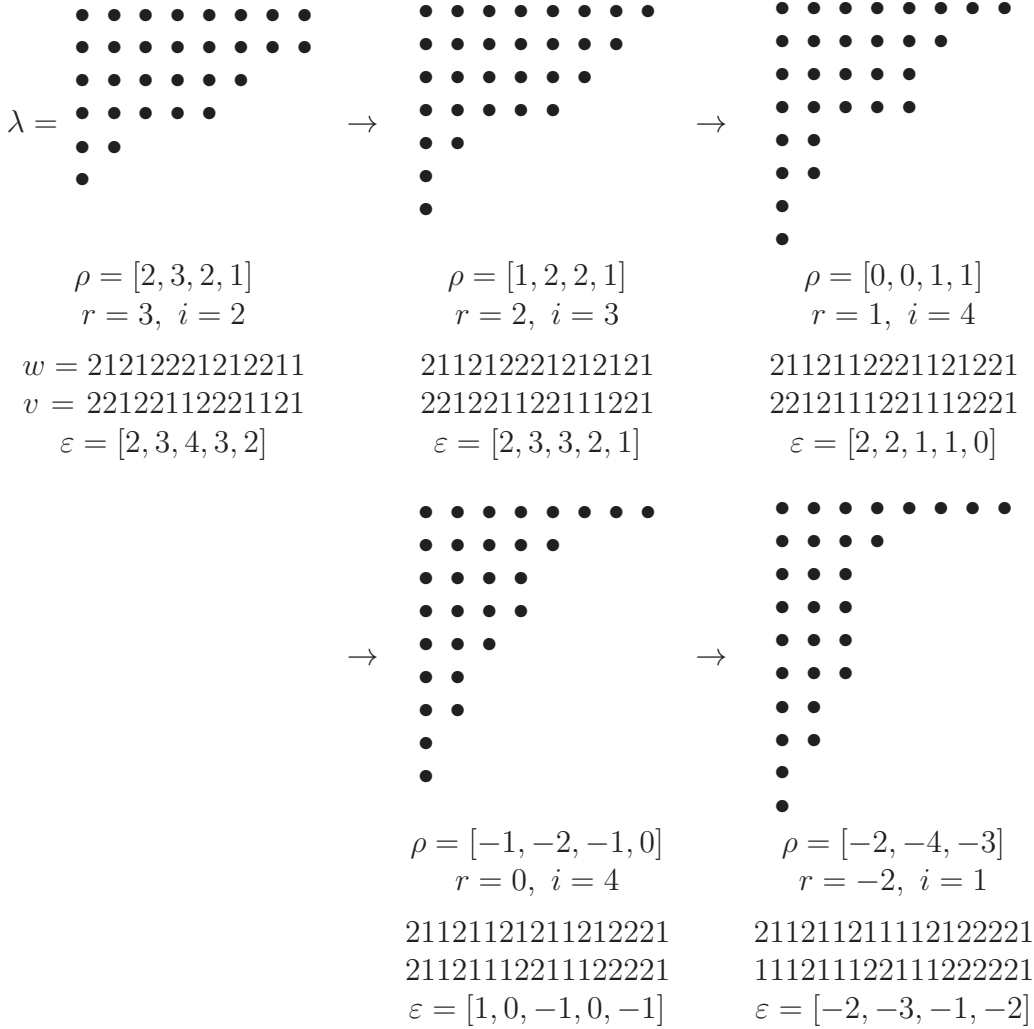
\thicklines
\setlength{\unitlength}{2pt}
$$
\barr{ccccc}
\la=
\barr{cccccccc}
\x & \x & \x & \x & \x & \x & \x & \x \\[-2pt]
\x & \x & \x & \x & \x & \x & \x & \x \\[-2pt]
\x & \x & \x & \x & \x & \x \\[-2pt]
\x & \x & \x & \x & \x \\[-2pt]
\x & \x \\[-2pt]
\x\\[20pt] 
\earr
&
\ra
&
\barr{cccccccc}
\x & \x & \x & \x & \x & \x & \x & \x \\[-2pt]
\x & \x & \x & \x & \x & \x & \x  \\[-2pt]
\x & \x & \x & \x & \x & \x \\[-2pt]
\x & \x & \x & \x & \x \\[-2pt]
\x & \x \\[-2pt]
\x \\[-2pt]
\x\\[10pt] 
\earr
&
\ra
&
\barr{cccccccc}
\x & \x & \x & \x & \x & \x & \x & \x \\[-2pt]
\x & \x & \x & \x & \x & \x \\[-2pt]
\x & \x & \x & \x & \x \\[-2pt]
\x & \x & \x & \x & \x \\[-2pt]
\x & \x \\[-2pt]
\x & \x \\[-2pt]
\x \\[-2pt]
\x 
\earr
\\
\rho=[2,3,2,1]
&&
\rho=[1,2,2,1]
&&
\rho=[0,0,1,1]
\\
 r=3,\ i=2
&&
r=2,\ i=3
&&
r=1,\ i=4
\\[5pt]
w=21212221212211
&&
211212221212121
&&
2112112221121221
\\
v\hs{1pt}=\hs{1pt}22122112221121
&&
221221122111221
&&
2212111221112221
\\
\vep=[2,3,4,3,2]
&&
\vep=[2,3,3,2,1]
&&
\vep=[2,2,1,1,0]
\\[10pt]
&
\ra
&
\barr{cccccccc}
\x & \x & \x & \x & \x & \x & \x & \x \\[-2pt]
\x & \x & \x & \x & \x \\[-2pt]
\x & \x & \x & \x \\[-2pt]
\x & \x & \x & \x \\[-2pt]
\x & \x & \x \\[-2pt]
\x & \x \\[-2pt]
\x & \x \\[-2pt]
\x \\[-2pt]
\x \\[10pt]
\earr
&
\ra
&
\barr{cccccccc}
\x & \x & \x & \x & \x & \x & \x & \x \\[-2pt]
\x & \x & \x & \x \\[-2pt]
\x & \x & \x \\[-2pt]
\x & \x & \x \\[-2pt]
\x & \x & \x \\[-2pt]
\x & \x & \x \\[-2pt]
\x & \x \\[-2pt]
\x & \x \\[-2pt]
\x \\[-2pt]
\x 
\earr
\\
&&
\rho=[-1,-2,-1,0]
&&
\rho=[-2,-4,-3]
\\
&&
r=0,\ i=4
&&
r=-2,\ i=1
\\[5pt]
&&
21121121211212221
&&
211211211112122221
\\
&&
21121112211122221
&&
111211122111222221
\\
&&
\vep=[1,0,-1,0,-1]
&&
\vep=[-2,-3,-1,-2]
\earr
$$

Let $\ka$ denote the map obtained by going once through the loop at CSV2.  Now for an arbitrary partition $\la$, $\ka(\la)$ may not be well defined since it may not be possible to find a part of size $i(\la)$ in $\la'$.  But it is shown in~\cite{csv:bpa} that if one starts with $\la\in D_0$ then such a part must exist while $r(\la)\ge0$.  One also needs to worry about termination of the algorithm.  But in~\cite{csv:bpa} it is shown that while $r(\la)\geq0$ we have $r(\ka(\la))\le r(\la)-1$ with equality if $r(\la)>0$.  So 
\beq
\label{CSV}
CSV(\la)=\ka^{r+1}(\la) 
\eeq
where $r=r(\la)$.

In their seminal paper~\cite{gk:svs}, Green and Kleitman gave a symmetric chain decomposition of the Boolean algebra $\cB_n$ of all subsets of $\{1,\ldots,n\}$ as follows.  Represent an element $S\in\cB_n$ as a word $w=w_1\ldots w_n\in\{1,2\}^*$ where $w_i=2$ if and only if $i\in S$.  Pair ones and twos in $w$ as described near the end of Section~\ref{bs}.  Note that the unpaired twos must precede the unpaired ones.  Define a map $\ga$ on words with at least one unpaired two by 
$$
\ga(w)=\mbox{$w$ with the rightmost unpaired two changed to a one.} 
$$
Note that, by the choice of the changed two, $\ga(w)$ has the same pairs as $w$.  The chains in the decomposition are all those of the form $w,\ga(w),\ga^2(w),\ldots,\ga^t(w)$ where the unpaired elements of $w$ are all twos and $t$ is the number of such elements (so the unpaired elements of $\ga^t(w)$ are all ones).

Define a map $GK:W_{121}\ra B_{21}$ as follows.  If $v=x121\in W_{121}$ has $t$ unpaired twos (note that they must all be in $x$), then let
\beq
\label{GK}
GK(v)=\ga^t(x) 12^{t+1} 1.
\eeq
\bth
The map $GK:W_{121}\ra B_{21}$ is well defined and bijective.
\eth
\prf
For $GK$ to be well defined we must check $GK(v)\in B_{21}$ for $v\in W_{121}$.  It is clear from~\ree{GK} that $v\in\{1,2\}^*21$.  Also, $GK(v)$ has no unmatched twos since $\ga^t(x)$ has replaced the $t$ unmatched twos in $x$ with ones, and so they can be used (along with the penultimate one) to match the twos in $2^{t+1}$.  But having no unmatched twos is equivalent to being a ballot sequence and so $GK(v)\in B_{21}$.

To prove bijectivity, we construct an inverse map.  Given $w\in B_{21}$ we write $w=y12^{t+1}1$ for some $t\ge0$.  Since $w$ is a ballot sequence, $t$ of the twos in the last $2$-run must be matched with ones in $y$ which are unmatched in $y$ itself.  So $x=\ga^{-t}(y)$ is well defined and we can set $v=x121$.  It is straightforward to check that $GK(v)=w$ and that this construction provides an inverse to $GK$.
\hqedm

We now have all the tools in hand to demonstrate the relationship between $CSV$ and $GK$.
\bth
We have
\beq
\label{CSV-GK}
CSV=\phi\circ GK\circ \phi^{-1}.
\eeq
\eth
\prf
Consider $\la\in D_0$ and $v=\phi^{-1}(\la)$.  If $\la\in R_{\le-1}$ as well then $CSV(\la)=\la$.  We also have $v\in W_{121}\cap B_{21}$.
It follows that $v=x121$ where $x$ is a ballot sequence and so has no unpaired twos.  This implies that
$GK(v)=v$ and both sides of~\ree{CSV-GK} agree, as desired.  

From now on we can assume $\la\not\in R_{\le-1}$ and thus $r=r(\la)\ge0$.  Appealing to Lemma~\ref{e>=r+1}, we have $e(v)\ge1$.  In this case  we will show  that if  $t$ is the number of unpaired ones in $v$ then $r+1=t$.  Note that for any $v$, if $e_i(v)$ is a positive left-right maximum in the sequence $e_1(v), e_2(v),\ldots$ then the $i$th $2$-run in $v$ has exactly $e_i(v)-e_j(v)$ unpaired twos, where $e_j(v)$ is the positive left-right maximum just before $e_i(v)$ (or $0$ if there is no such previous value).  Since $e(v)>0$ we know that a positive left-right maximum exists and thus $e(v)=t$, the number of unpaired twos.  We also have that $v\in W_{121}$ implying that  $e_d(v)=e_{d-1}(v)-1$ and so $e_d(v)<e(v)$.  Applying Lemma~\ref{e>=r+1} again shows that $r+1=e(v)=t$ as claimed.

We now consider a variant $\gab$ of $\ga$ as follows.  Let $u=x12^m1$ where $x$ has at least one unpaired two and $m\ge1$. Define
$$
\gab(u)=\ga(x)12^{m+1}1.
$$
By the definitions~\ree{CSV} and~\ree{GK}, we have $GK(v)=\gab^t(v)$ and $CSV(\la)=\ka^{r+1}(\la)$ where $t=r+1$.  So the theorem will follow if we can show that $\ka=\phi\circ\gab\circ\phi^{-1}$ or, equivalently, that $\ka(\phi(y))=\phi(\gab(y))$ whenever $y=\gab^j(v)$ for some $j<t$.

To set notation, let $\mu=\la(\phi(y))$, $z=\gab(y)$, and $\nu=\la(\phi(z))$.  We can complete the proof by showing that $\ka(\mu)=\nu$.  Let $y$ have the form~\ree{des}.  Since $y=\gab^j(v)$ where $v\in W_{121}$ we see that $n_d=0$, $m_d=1$, and $n_{d-1}=j+1$.  We must also locate the rightmost unpaired two in $y$.  As in the CSV algorithm, let $i=i(\mu)$ denote the maximum index such that $r_i(\la)=r(\la)$.  Using equation~\ree{e_i(v)} and the fact that $e_d(y)<e(y)$ since $y$ ends with a single one, we see that $d-i$ is the minimum index such that $e_{d-i}(y)=e(y)$.  So, by the description of the location of the unpaired twos in the second paragraph of this proof,  the rightmost such must lie in the $(d-i)$th run.  Combining all the information gathered thus far and  using the definition of $\gab$ permits us to write
\bal
y&=p 2^{n_{d-i}} 1^{m_{d-i+1}} q 1 2^{j+1} 1,\\
z&=p 2^{n_{d-i}-1} 1^{m_{d-i+1}+1} q 1 2^{j+2} 1
\end{align*}
for some words $p,q$.

We now wish to apply $\phi$ using Proposition~\ref{phi}.  To do so, we must worry about whether $n_{d-i}-1=0$ in $z$.  But in that case $n_{d-i}=1$ and the only way a run consisting of a single two could have that two unpaired is if it is the initial run.  Since this is also the rightmost unpaired two, it must be the only unpaired two and so $j=t-1$.  The proof in this special case is similar to the one for $j<t-1$ which we will present, leaving the last value of $j$ for the reader.

For $j<t-1$ we can write
\bal
\phi(y)&=\ol{p} 1^{m_{d-i+1}-1} \ol{q}  2^{n_{d-i}-1}\ol{s} 1 2^{j} 1,\\
\phi(z)&=\ol{p} 1^{m_{d-i+1}} \ol{q}  2^{n_{d-i}-2}\ol{s} 1 2^{j+1} 1
\end{align*}
for certain $\ol{p}, \ol{q}, \ol{s}$.  So to pass from $\phi(y)$ to $\phi(z)$ we have done the following.  We have added a one to the run $1^{m_{d-i+1}-1}$ in $\phi(y)$.  By Proposition~\ref{phi}, there are exactly $i-1$ twos before this run and so this corresponds to adding a part of size $i-1$ to $\mu$ when forming $\nu$.  Similarly, subtracting a two from the run $2^{n_{d-i}-1}$ corresponds to deleting a part of size $i$ from $\mu'$.  And a two is added to the final $2$-run, which corresponds to increasing the size of $\mu_1$ by one.  In summary, we have $\ka(\mu)=\nu$ and are done. 
\hqed

\section{The past, the future, and open problems}
\label{pfo}

\subsection{Lucanomials}

We will describe the genesis of the notion of a Mahonian pair.  This is not only for historical reasons, but also because the original problem which lead us to this definition is still unsolved.

Let $s,t$ be variables and define a sequence of polynomials $\{n\}$ in $s$ and $t$ by $\{0\}=0$, $\{1\}=1$, and 
\beq
\label{luc}
\{n\} = s\{n-1\} + t\{n-2\}
\eeq
for $n\ge 2$.   When $s,t$ are integers, the arithmetic properties of this sequence were studied by Lucas~\cite{luc:tfn1,luc:tfn2,luc:tfn3}.  Given integers $0\le k\le n$, define the corresponding {\it lucanomial coefficient\/}
$$
\bin{n}{k}=\frac{\{n\}!}{\{k\}!\{n-k\}!}.
$$
Note that one obtains the fibonomial coefficients or $q$-binomial coefficients by specializing to $s=t=1$ or to $s=[2]$ and $t=-q$, respectively.
It is not hard to show that these rational functions are actually polynomials in $s,t$ with coefficients in $\bbN$.  In~\cite{ss:cib}, we were able to give a simple combinatorial interpretation of the lucanomial coefficients using tilings of partitions contained in a $k\times(n-k)$ rectangle.  (Earlier, more complicated, interpretations were given by Gessel-Viennot~\cite{gv:bdp} and Benjamin-Plott~\cite{bp:caf,bp:caferr}.)

Upon hearing a talk on this subject at the 2010 Mathfest, Lou Shapiro asked the following natural question.  Define an $s,t$-analogue of the Catalan numbers by
$$
C_{\{n\}}=\frac{1}{\{n+1\}}\bin{2n}{n}.
$$
\begin{question}[Shapiro]
Is $C_{\{n\}}\in\bbN[s,t]$?  If so, is there a simple combinatorial interpretation?
\end{question}
Ekhad~\cite{ekh:ssl} has pointed out that the answer to the first question is ``yes"  since one can show that
$$
C_{\{n\}}= \bin{2n-1}{n-1}+t\bin{2n-1}{n-2}.
$$
However, a nice combinatorial interpretation is elusive.  While trying to find such an interpretation, we looked at statistics on ballot sequences (being one of the most common objects associated with $C_n$) and partitions with all ranks positive sitting inside an $n\times n$ rectangle (since Savage had done previous work with such partitions).  It was noted that $\maj$ on the former was equidistributed with area on the latter, and the concept of a Mahonian pair was born.

\subsection{Larger alphabets}

It would be very interesting to obtain results about Mahonian pairs using sets outside of $\{1,2\}^*$.  One natural  place to look would be at sets of permutations determined by pattern avoidance.  If $\pi\in\fS_k$ then we say $\si\in\fS_n$ contains $\pi$ as a {\it pattern\/} if there is a subpermutation of $\si$ order isomorphic to $\pi$.  If $\si$ does not contain $\pi$ as a pattern then we say it {\it avoids\/} $\pi$.  We let $\Av_n(\pi)$ denote the set of such patterns in $\fS_n$.   Since $|\Av_n(\pi)|=C_n$ for any $\pi\in\fS_3$, the hope was that there would be some connection with the Mahonian pairs associated with ballot sequences.  Unfortunately, looking at the distributions for small $n$ showed no possible pairs among the $\Av_n(\pi)$ where $\pi\in\fS_3$.

However, another phenomenon manifested itself.  There were pairs $\pi,\si$ such that the distribution of $\inv$ over both $\Av_n(\pi)$ and $\Av_n(\si)$ were the same, and similarly for $\maj$.   This lead us to the following refinement of the traditional notion of Wilf equivalence.  (We say $\pi$ and $\si$ are {\it Wilf equivalent\/} if $|\Av_n(\pi)|=|\Av_n(\si)|$ for all $n\ge0$.)  Let $\sta$ be any statistic on permutations.  Call $\pi$ and $\si$ {\it st-Wilf equivalent\/} if $\sta$ is equidistributed over $\Av_n(\pi)$
and $\Av_n(\si)$.  So st-Wilf equivalence implies Wilf-equivalence, but not conversely.  

In very recent work, Dokos, Dwyer, Johnson, Sagan, and Selsor~\cite{ddjss:pps} studied this concept for the $\inv$ and $\maj$ statistics.  They also considered, for $\Pi\sbe\fS_3$, the sets $\Av_n(\Pi)=\cap_{\pi\in\Pi}\Av_n(\pi)$.  For such sets when $|\Pi|\ge2$, they found quite a number of Mahonian pairs.  Here is a sample.
\bth[\cite{ddjss:pps}]
Let $S=\Av_n(\Pi)$ where
$$
\Pi\in\large\{\{132,213\}, \{132, 312\}, \{213, 231\}, \{231,312\}\large\}
$$
and let $T=\Av_n(\Pi')$ where
$$
\Pi'\in\large\{\{132,231\}, \{132, 312\}, \{213, 231\}, \{213,312\}\large\}.
$$
Then $(S,T)$ form a Mahonian pair.\hqedm
\eth

\subsection{Eulerian pairs}

There is another famous pair of equidistributed statistics.  An {\it excedance\/} in a permutation $w=a_1 \ldots a_n$  is an index $i$ such that $a_i>i$.  Let $\exc w$ be the number of excedences.  It is well known that $\des$ and $\exc$ are equidistributed over $\fS_n$ and any statistic with this distribution is called {\it Eulerian\/}.  It is easy to extend the definition of excedance to any $w\in\bbP^*$.  Let $x=b_1\ldots b_n$ be the weakly increasing rearrangement of $w$.  Then an excedance of $w$ is an index $i$ with $a_i>b_i$.  The next definition should come as no surprise.   The subsets $S,T\sbe\bbP^*$ form an {\it Eulerian pair\/}, $(S,T)$, if there is a bijection $\al:S\ra T$ such that 
$$
\exc v = \des \al(v)
$$
for all $v\in S$.  The study of Eulerian pairs could be every bit as rich as that of their Mahonian cousins and we will be investigating their properties.

{\it Acknowledgment.\/}  We would like to thank Dominique Foata, Ira Gessel,  Drew Sills, and Michelle Wachs for useful discussions and references.

\end{document}